\documentstyle[12pt]{article}
\def\Tilde{\char126\relax}

\begin{document}
\begin{center}\Large{New Lower Bounds for Some Multicolored Ramsey Numbers}
\end{center}
\vskip 10pt
\begin{center}Aaron Robertson\footnote{webpage: 
www.math.temple.edu/\Tilde aaron\\ \indent This paper is part of the
author's Ph.D. thesis under the direction of Doron Zeilberger.\\
\indent This paper was supported in part by the NSF under the
PI-ship of Doron Zeilberger.}\\
Department of Mathematics, Temple University\\ Philadelphia, PA
19122\\
email:  aaron@math.temple.edu\\
\vskip 5pt
Submitted:  October 26, 1998;  Accepted November 15, 1998\\
Classification:  05D10, 05D05
\end{center}
\vskip 30pt

\begin{abstract}
\noindent
In this article we use two different methods to find
new lower bounds for some multicolored Ramsey numbers.
In the first part we use the finite field method used by
Greenwood and Gleason [GG] to show that
$R(5,5,5) \geq 242$ and $R(6,6,6) \geq 692$. 
In the second part
we extend Fan Chung's result in [C] to show that, 
$$
R(3,3,3,k_1,k_2,\dots,k_r) \geq 3 R(3,3,k_1,k_2,\dots,k_r) +
R(k_1,k_2,\dots,k_r) - 3
$$ 
holds for any natural number $r$ and for any $k_i\geq 3$,
$i=1,2,\dots r$. This general result, along with known results,
imply the following nontrivial bounds:
$R(3,3,3,4) \geq 91$, $R(3,3,3,5) \geq 137$, 
$R(3,3,3,6) \geq 165$, $R(3,3,3,7) \geq 220$,
$R(3,3,3,9) \geq 336$, and $R(3,3,3,11) \geq 422$.
\end{abstract}
\vskip 20pt

\begin{center}
\textbf{Introduction}
\end{center}

This paper is presented in two part which can be read independently of
each other.  Part one uses finite fields, and Part two extends an
argument by Fan Chung.
\vskip 10pt
Part one of this article is accompanied by the Maple package {\tt RES},
available for
download at the author's website.
\vskip 10pt
Recall that $N=R(k_1,k_2,\dots,k_r)$ is the minimal integer with
the following property:
\vskip 10pt
\hskip 20pt {\bf Ramsey Property}: {\it If we $r$-color the edges of
the complete graph on $N$ vertices,
\vskip 2pt
\hskip 20pt $K_N$, then there {\bf must}
exist a $j$, $1 \leq j \leq r$, such that a monochromatic $j$-colored
\vskip 2pt
\hskip 20pt
complete graph on $k_j$ vertices is a subgraph of the $r$-colored
$K_N$.}

\newpage

\begin{center}
\textbf{Part One:  The Finite Field Method}
\end{center}

In this first part we 
add two more lower bounds to Radziszowski's Dynamic Survey [R] on the
subject.
We show, by using the finite field technique in [GG], that 
$R(5,5,5) \geq 242$ and $R(6,6,6) \geq 692$.  The previous best lower bound
for
$R(5,5,5)$ was $169$ given by Song [S], who more generally shows that
$R(\underbrace{5,5, \dots, 5}_{r \ times}) \geq 4(6.48)^{r-1} + 1$ holds
for all $r$.  For $R(6,6,6)$ there was no established nontrivial lower bound.
\vskip 10pt
Consider the number $R(5,5,5)$.  To find a lower bound, $L$, 
we are searching for a three coloring of $K_L$ which avoids a 
monochromatic $K_5$.  We use an argument of Greenwood and Gleason, which
is reproduced here for the sake of completeness.
\vskip 10pt
Let $L$ be prime and consider the field of $L$ elements,
numbered from $0$ to $L-1$.
Associate each field element with a vertex of $K_L$.  We require
that $3$ divides $L-1$.  Now consider the cubic 
residues of the multiplicative group
${\mathbf{Z}}_L^* = {\mathbf{Z}}_L \setminus \{0\}$, which form a coset of 
 ${\mathbf{Z}}_L^*$. 
Note that since $3$ divides $L-1$, there must be $2$ other cosets.  
\vskip 10pt
Let $i$ and $j$ be two vertices of $K_L$.  Color the edges of $K_L$ as
follows:  If $j-i$ is a cubic residue color the edge connecting 
$i$ and $j$ red, if it is in the second coset, color the edge blue, and 
if it is the third coset, color the edge green.  (Note that 
the order of differencing is immaterial since $-1$ is a cubic residue.) 
\vskip 10pt
Now suppose that a monochromatic $K_5$ exists in this coloring.  Without
loss of generality we may call the five vertices $0$, $a$, $b$, $c$, and 
$d$, with $0<a<b<c<d$.
Then the set of edges, $E= \{a,b,c,d,b-a,c-a,d-a,c-b,d-b,d-c\}$, must
be a subset of one of the cosets.  Since $a \neq 0$,
multiplication by $a^{-1}$ is allowed.  Set $B=ba^{-1}$, $C=ca^{-1}$, 
and $D=da^{-1}$.  Then the set $a^{-1}E= \{1,B,C,D,B-1,C-1,D-1,C-B,D-B,D-C \}$
must be a subset of the cubic residues.  Hence if we find an $L$ for which
there does not exist $B,C,$ and $D$ such that $a^{-1}E$ is a subset of the
cubic residues, then we can conclude that $R(5,5,5) > L$.
Of course, this argument holds for $
R(\underbrace{t,t, \dots, t}_{k \ times})$ for any $k$, and any $t$.
\vskip 20pt
\begin{center}
\textbf{Using {\tt RES}}
\end{center}

We only acheived results when we restricted our search to fields of
prime order (although any finite field can be explored using {\tt RES}
(or at least easily modified to do so)).  Since we are considering the
number $R(5,5,5)$, reject
 any prime, $q$, for which $3$ does not divide $q-1$.  This can be
accomplished automaticly by using the procedure {\tt pryme}. 
By using the procedure {\tt res}
we produce all of the cubic residues of ${\mathbf{Z}}_p^*$, 
for a given prime, $p$.  We then
use the procedure {\tt siv} to discard any residue, $R$, for which 
$R-1$ is not a residue.  We now have a much more manageable list to search.
We then call the procedure {\tt diffcheck} to choose all possible $3$-sets
(for $B$, $C$, and $D$)
and check whether or not the differences between any two elements are
\emph{all}
cubic residues.  If such a $3$-set exists, {\tt diffcheck} will output 
the first $3$-set it finds.  However, in the event that no such $3$-set
exists, {\tt diffcheck}
will output $1$. 
\vskip 5pt
{\tt RES} can also be used to search finite fields whose order
is not prime.  For example, to verify that the field on $2^4$ elements,
avoids a monochromatic triangle by using cubic residues (this fact was
proven in [GG]), type {\tt GalField3(2,4,3)}. 
\vskip 10pt
By using {\tt RES} we were able to find the following lower bounds:
$R(5,5,5) \geq 242$ and $R(6,6,6) \geq 692$.  These are obtained by the
following colorings: (Since $-1$ is a cubic residue it
suffices to list only entries up to $120$ for $R(5,5,5)$ and
$345$ for $R(6,6,6)$.)
\vskip 10pt
$R(5,5,5) > 241$:
\vskip 10pt
{\tt {\bf Color 1:} 1, 5, 6, 8, 17, 21, 23, 25, 26, 27, 28, 30, 33,
36, 40, 41, 43, 44, 47, 48, 57, 61, 64, 73, 76, 79, 85, 87, 91, 93,
98, 101, 102, 103, 105, 106, 111, 115, 116, 117}
\vskip 10pt
{\tt {\bf Color 2:}  2, 7, 9, 10, 11, 12, 16, 19,
29, 31, 34, 35, 37, 39, 42, 45, 46, 50, 52, 54, 55, 56, 59, 60,
66, 67, 71, 72, 80, 82, 83, 86, 88, 89, 94, 95, 96, 113, 114, 119}
\vskip 10pt
{\tt {\bf Color 3:}  3, 4, 13, 14, 15, 18, 20, 22, 24, 32, 38, 49,
51, 53, 58, 62, 63, 65, 68, 69, 70, 74, 75, 77, 78, 81, 84, 90, 92,
97, 99, 100, 104, 107, 108, 109, 110, 112, 118, 120}
\vskip 10pt
$R(6,6,6) > 691$:
\vskip 10pt
{\tt {\bf Color 1:} 1, 2, 4, 5, 8, 10, 16, 19, 20, 21, 25, 27, 31,
32, 33, 38, 39, 40, 42, 50, 51, 54, 62, 64, 66, 67, 69, 71, 73, 76,
78, 80, 83, 84, 87, 89, 95, 100, 102, 105, 107, 108, 109, 123, 124,
125, 128, 132, 134, 135, 138, 139, 142, 146, 149, 151, 152, 155,
156, 160, 163, 165, 166, 168, 173, 174, 178, 179, 181, 190, 191,
195, 199, 200, 204, 210, 214, 216, 218, 246, 248, 250, 255, 256,
259, 263, 264, 268, 270, 271, 276, 278, 283, 284, 291, 292, 293,
298, 301, 302, 304, 309, 310, 311, 312, 320, 326, 329, 330, 332,
333, 335, 336, 343, 345}
\vskip 10pt
{\tt {\bf Color 2:}  7, 9, 11, 13, 14, 17, 18, 22, 23, 26, 28, 29,
34, 35, 36, 41, 44, 45, 46, 52, 55, 56, 58, 65, 68, 70, 72, 82, 85,
88, 90, 92, 97, 103, 104, 110, 111, 112, 115, 116, 127, 129, 130,
131, 133, 136, 140, 141, 144, 145, 147, 159, 164, 167, 170, 171,
175, 176, 177, 180, 183, 184, 189, 194, 197, 205, 206, 208, 209,
217, 220, 222, 224, 225, 227, 229, 230, 231, 232, 233, 237, 241,
243, 247, 251, 254, 257, 258, 260, 262, 266, 272, 273, 275, 279,
280, 281, 282, 288, 290, 294, 297, 303, 313, 318, 323, 325, 328,
331, 334, 337, 339, 340, 341, 342}
\vskip 10pt
{\tt {\bf Color 3:}  3, 6, 12, 15, 24, 30, 37, 43, 47, 48, 49, 53,
57, 59, 60, 61, 63, 74, 75, 77, 79, 81, 86, 91, 93, 94, 96, 98, 99,
101, 106, 113, 114, 117, 118, 119, 120, 121, 122, 126, 137, 143,
148, 150, 153, 154, 157, 158, 161, 162, 169, 172, 182, 185, 186,
187, 188, 192, 193, 196, 198, 201, 202, 203, 207, 211, 212, 213,
215, 219, 221, 223, 226, 228, 234, 235, 236, 238, 239, 240, 242,
244, 245, 249, 252, 253, 261, 265, 267, 269, 274, 277, 285, 286,
287, 289, 295, 296, 299, 300, 305, 306, 307, 308, 314, 315, 316,
317, 319, 321, 322, 324, 327, 338, 344}
\vskip 40pt
\begin{center}
{\bf Part Two:  On the Ramsey Numbers $R(3,3,3,k_1,k_2,\dots,k_r)$}
\end{center}
\vskip 10pt

Let $N=R(k_1,k_2,\dots,k_r)$.  The Ramsey Property
implies that there must exist a graph on $N-1$ vertices
which avoids the Ramsey Property.  Using such a graph, along with
the construction in [C], we will prove that, for any natural number
$r$ and for any $k_i \geq 3$, $i=1,2,\dots r$, 
$$
R(3,3,3,k_1,k_2,\dots,k_r) \geq 3 R(3,3,k_1,k_2,\dots,k_r) +
R(k_1,k_2,\dots,k_r) - 3.
$$ 
\vskip 20pt

\begin{center}
\textbf{The Construction}
\end{center}

Fix $r \geq 1$, and $k_i \geq 3$ for $i=1,2,\dots r$.  Let
$M=R(3,3,k_1,k_2,\dots k_r)$.  Then there exists a graph, $G$, on
$M-1$ vertices which avoids the Ramsey Property.  Call the incidence
matrix of this graph $T_{r+2}=T_{r+2}(x_0,x_1,x_2,\dots, x_{r+2})$, where
$x_0$ are the diagonal entries only, and the $x_i$, for $i=1,2, \dots
r+2$, are the $r+2$ colors.  By definition of $G$, there are no
$x_1$-colored nor $x_2$-colored triangles, and no $x_{i+2}$-colored
$K_{k_i}$, for $i=1,2, \dots r$.
\vskip 10pt

Now consider the following slightly modified construction from [C]:
\vskip 10pt
$$
\begin{array}{ccccc}
&A\\
&D&B\\
T_{r+3}(0,1,2,\dots,r+3)=&E&F&C\\
&1,\dots,1&2,\dots,2&3,\dots,3\\
&\vdots&\vdots&\vdots&G\\
&1,\dots,1&2,\dots,2&3,\dots,3\\
\end{array}
$$
\vskip 10pt
\noindent
the incidence matrix of a graph $H$ on $3M+R(k_1,k_2,\dots, k_r)-4$
vertices, where 
$$
\begin{array}{c}
A=T_{r+2}(0,2,3,4,5,\dots,r+3)\\
B=T_{r+2}(0,3,1,4,5,\dots,r+3)\\
C=T_{r+2}(0,1,2,4,5,\dots,r+3)\\
D=T_{r+2}(3,2,1,4,5,\dots,r+3)\\
E=T_{r+2}(2,1,3,4,5,\dots,r+3)\\
F=T_{r+2}(1,3,2,4,5,\dots,r+3)\\
\end{array}
$$
\noindent
and $G$ is any matrix on $R(k_1,k_2,\dots,k_r)-1$ vertices in the colors
$4$ through $r+3$ which avoids the
Ramsey Property.  
\vskip 5pt
Using Fan Chung's result [C] we see that the graph $H$ avoids
$1$-colored, $2$-colored, and 
$3$-colored triangles.  We now argue that no $(j+3)$-colored 
$K_{k_j}$ exists in $H$ for $j=1,2, \dots, r$:  Assume there exists
a $J$-colored $K_{k_J}$ in H, for some $J$ between $4$ and $r+3$. 
Then there must exist ${k_J \choose 2}$ entries in the $3(M-1) \times 3(M-1)$
upper left submatrix of $T_{r+3}$, all
of value $J$, which form the edges of a complete graph on $k_J$ 
vertices.  However, 
by the construction of $T_{r+3}$ we see that {\it all} of these
entries can be taken modulo $M-1$, since the entries of value $J$ in each
block are in exactly the same places as in the upper left block, $A$. 
We further note that if
$(s,t)$ and $(u,v)$ are two of the entries in question, then 
$(s,t) \not\equiv (u,v) \ (mod \ M-1)$ (componentwise).  Without loss
of generality we may assume $s < u$.  
If we had $s \equiv u \ (mod \ M-1)$, then since $(u,s)$ must also have
the same value as $(u,v)$, we would have $(u,s) \equiv (u,u) \ (mod \ M-1)$.
This implies that the entry $J$ is on the diagonal of 
$T_{r+2}(0,2,3,4, \dots, r+3)$, a contradiction. Hence, if we have a 
$J$-colored $K_{k_J}$ in the upper left submatrix of  
$T_{r+3}$ (of size $3(M-1) \times 3(M-1)$), then there must be a $J$-colored
$K_{k_J}$ in $T_{r+2}(0,2,3,4, \dots, r+3)$, contradicting the
definition of $T_{r+2}$.
\vskip 5pt
\noindent
{\it Remark}:  Up to the renaming of colors and vertices, the above
permutation configuration of colors which defines $T_{r+3}$ is the
{\it only} configuration which will avoid monochromatic triangles. 
\vskip 20pt

\begin{center}
\textbf{Harvesting Some Lower Bounds for Ramsey Numbers}
\end{center}

It is amazing that the result of this section has not be
observed for the $25$ years since [C] was published.  
Using this observation we will give $6$ new lower bounds.
Currently in Radziszowski's Survey [R], we have that
$R(3,3,3,4) \geq 87$, due to Exoo [E1].
By applying the result of this section to
$R(3,3,3,4)$ and using the fact that $R(3,3,4) \geq 30$ [K], we get
the new lower bound:  $R(3,3,3,4) \geq 91$.
Using the bound $R(3,3,5) \geq 45$ [E2,KLR], we get the bound 
$R(3,3,3,5) \geq 137$.  Finally, using $R(3,3,6) \geq 54$, 
$R(3,3,7) \geq 72$, $R(3,3,9) \geq 110$, and 
$R(3,3,11) \geq 138$ all from [SLZL], we get the lower bounds
$R(3,3,3,6) \geq 165$, $R(3,3,3,7) \geq 220$, $R(3,3,3,9) \geq 336$, 
and $R(3,3,3,11) \geq 422$.
\vskip 10pt
\noindent
\textbf{Acknowledgments} 
\vskip 10pt
I would like to thank my advisor, Doron Zeilberger, for his
useful comments and insight regarding this paper.  More importantly, however, 
I would like to thank him for sparking my interest in 
combinatorics and for sharing his mathematical philosophies.  The
mathematical community could benefit greatly from more mathematicians
like him. 

I would also like to thank Brendan McKay for his assistance and 
the referee for calling my attention to [SLZL].
\vskip 20pt
\noindent
\textbf{References}
\vskip 10pt
\noindent
[C] F. Chung, {\it On the Ramsey Numbers $N(3,3, \dots, 3;2)$},
Discrete Mathematics, {\bf 5}, 1973,\\
\indent
317-321.
\vskip 5pt
\noindent
[E1] G. Exoo, {\it Some New Ramsey Colorings}, Electronic Journal of
Combinatorics, R29, {\bf 5},\\
\indent
1998, 5pp.
\vskip 5pt
\noindent
[E2] G. Exoo, {\it Constructing Ramsey Graphs with a Computer},
Congressus Numerantium, {\bf 59},\\
\indent
1987, 31-36.
\vskip 5pt
\noindent
[GG] R. Greenwood and A. Gleason, {\it Combinatorial Relations and 
Chromatic}\\ 
\indent
{\it Graphs}, Canadian Journal of Mathematics, \textbf{7}, 1955, 1-7.
\vskip 5pt
\noindent
[K] J. Kalbfleisch, {\it Chromatic Graphs and Ramsey's Theorem}, 
(Ph.D. Thesis), University of\\
\indent
Waterloo, January 1966.
\newpage
\noindent
[KLR] D.L. Kreher, Li Wei, and S. Radziszowski, {\it Lower Bounds for
Multi-Colored Ramsey}\\ 
\indent
{\it Numbers from Group Orbits}, Journal of Combinatorial Mathematics 
and Combinatorial\\
\indent
Computing, {\bf 4}, 1988, 87-95.
\vskip 5pt
\noindent
[R] S. Radziszowski, {\it Small Ramsey Numbers}, Electronic 
Journal of Combinatorics, Dynamic\\ 
\indent
Survey {\bf DS1}, 1994, 28pp.
\vskip 5pt
\noindent
[S] Song En Min, {\it New Lower Bound Formulas for the Ramsey Numbers
$N(k,k, \dots, k;2)$} (in\\ 
\indent
Chinese), Mathematica Applicata, {\bf 6}, 1993 suppl., 113-116.
\vskip 5pt
\noindent
[SLZL] Su Wenlong, Luo Haipeng, Zhang Zhengyou, and Li Guiqing, 
{\it New Lower Bounds of}\\     
\indent
{\it Fifteen Classical Ramsey Numbers}, to appear in Australasian
Journal of Combinatorics.
\end{document}